\documentclass{amsproc}

\newcommand{\Dem}{{\bf Proof: }}

\newtheorem{Pro}{Proposition}
\newtheorem{Cor}{Corollary}
\newtheorem{Lem}{Lemma}
\newtheorem{Teo}{Theorem}

\newcommand{\cqd}{\hfill$\Box$}

\newcommand{\Z}{{\mathbb Z}}
\newcommand{\R}{{\mathbb R}}
\newcommand{\C}{{\mathbb C}}

\newcommand{\ac}{{\mathfrak A}}
\newcommand{\Bc}{{\mathfrak B}}
\newcommand{\ec}{{\mathfrak E}}
\newcommand{\kc}{{\mathfrak K}}
\newcommand{\bc}{{\mathfrak B}}

\newcommand{\Sc}{{\mathfrak S}}

\title[Principal-symbol index for an algebra of pseudos]{The principal-symbol index map for an algebra of pseudodifferential operators}
\author{Severino T. Melo}
\address{Instituto de Matem\'atica e Estat\'istica, Universidade de S\~ao Paulo, 
Rua do Mat\~ao 1010,
05508-090 S\~ao Paulo, Brazil}
\email{toscano@ime.usp.br}
\dedicatory{This paper is dedicated to the memory of Heinz-Otto Cordes (1925-2018)}

\begin{document}

\begin{abstract}

  Let $B$ be a compact Riemannian manifold, let $\Omega$ denote the cylinder $\R\times B$, $\Delta_\Omega$ its Laplace operator
  and $\Lambda=(1-\Delta_\Omega)^{-1/2}$. Let $\mathfrak{A}$ denote the 
C*-algebra of bounded operators on $L^2(\R\times B)$ generated by all the classical  
pseudodifferential operators on $\R\times B$ of the form $L\Lambda^N$, $N$ a nonnegative integer and $L$ an
$N$-th order differential operator whose (local) coefficients approach $2\pi$-periodic functions at $+\infty$ and
$-\infty$. Let $\mathfrak{E}$ denote the kernel of the continuous extension of the
principal symbol to $\mathfrak{A}$.

The problem of computing the K-theory index map 
$\delta_1(K_1(\ac/\ec))\to K_0(\ec)\simeq\Z^2$ on an element of 
$K_1(\ac/\ec)$ is reduced to the problem of computing the Fredholm indices of two elliptic operators
on the compact manifold $S^1\times B$.

For $B=S^1$, Hess went further and proved in her thesis that  $K_0(\ac)\cong\Z^5$ and $K_1(\ac)\cong\Z^4$.

\end{abstract}

\maketitle

\section*{Introduction}

More than fifty years after the Atiyah-Singer index theorem, people are still looking for generalizations of their theorem for
classes of non-compact manifolds, for manifolds with boundary, corners or edges, or for manifolds with singularities in a more
general sense. There are, in fact, a great variety of different classes of problems, and of different approaches or tecnhiques that could be used.

One of the approaches is to consider C*algebras generated by zero-order pseudodiferential operators which contain an {\em order-reducing
device}. That means that there are linear isomorphisms between Sobolev spaces so that certain classes of differential operators
multiplied on the right by them fall into the C*algebra. In the case of a compact manifold,
the kernel of the continuous extension to the C*-algebra of the principal symbol map coincides with the ideal of compact operators,
and the Fredholm index can be regarded as the K-theory index map associated to the exact sequence (\ref{ases}) below. For a noncompact manifold,
one might argue that the relevant index is no longer the Fredholm index, but rather the K-theory index map taking values in $K_0$ of
the kernel of the principal symbol.

Monthubert and Nistor followed that track in \cite{MN1,MN2}, using algebras of pseudodifferential operators on groupoids which relate to (regular)
algebras of pseudodifferential operators on a noncompact manifold. In their case, the relevant compactification of the original manifold is itself
a manifold with corners. In this paper we consider an algebra of pseudodifferential operators on a cylinder $\R\times B$, where $B$ is compact, which includes
operators with periodic symbols. In our case the relevant compactification of the manifold is not even locally path connected.

The $K$-groups of our algebra, for the special cases of $B$ equals a point and of
$B$ equals $S^1$, have been computed in \cite{TO1} and \cite{Hess}, respectively.

The algebra considered here is an example of a {\em comparison algebra}, in the sense of Cordes,
who pioneered \cite{BC,CI,CO2,CH1,CH2} the use of Banach-algebra techniques in the study of the
spectral theory of partial differential operators and singular integral operators on noncompact manifolds.
Those algebras were defined in \cite{CO1} as C*-algebras of  bounded operators on $L^2(M)$
for a noncompact manifold $M$ equipped with a smooth positive measure and a positive second-order elliptic differential operator $H$. Typically,
$M$ is equipped with a riemannian metric and its induced measure, and $H$ is equal to the identity minus the Laplacian. One defines
$\Lambda=H^{-1/2}$, where $H$ also denotes the Friedrichs extension of the original operator defined on $C_c^\infty(M)$, which of course coincides with the unique
self-adjoint extension in case the operator is essentially self-adjoint. Given a class of smooth functions $\mathcal{A}^\sharp$ and a
class of first-order differential operators $\mathcal{D}^\sharp$ on $M$ satisfying certain conditions, a comparison algebra is defined as the C*-subalgebra of the
algebra of bounded operators on $L^2(M)$ generated by all operators of multiplication by $a\in \mathcal{A}^\sharp$ and by all
$D\Lambda$, $D\in \mathcal{A}^\sharp$. In his monograph \cite{CO1}, Cordes develops general techniques to describe the ideal structure of a comparison algebra
and how to find differential operators {\em within reach} of the algebra, meaning differential operators which, composed with
an isomorphism from $L^2$ to a suitable Sobolev space, belongs to the algebra.

In case the manifold is compact, there is only one possible choice of classes: $\mathcal{A}^\sharp=C^\infty(M)$ and $\mathcal{D}^\sharp$ equals
to the set of all first order differential operators with smooth coefficients. In this case, the resulting comparison algebra coincides \cite{T} with the
C*algebra generated by all zero-order classical pseudodifferential operators on $M$, which we will denote by $\Psi(M)$. The principal symbol of
operators, a priori a *-homomorphism defined on the subalgebra of classical pseudodifferential operators and
taking values in smooth functions on the co-sphere bundle $S^*M$,
extends to a surjective
C*-algebra homomorphism $\sigma_M:\Psi(M)\to C(S^*M)$ with kernel equal to the ideal of compact operators $\kc_M$. In other words, we have the
short exact sequence of C*-algebras
\begin{equation}
    \label{ases}
  0 \longrightarrow\kc_{M}\hookrightarrow\Psi(M)\stackrel{\sigma_{M}}{\longrightarrow} C(S^*(M))\longrightarrow 0.
  \end{equation}
In Cordes' approach, the map $\sigma_M$ arises as the Gelfand map for the commutative C*algebra $\Psi(M)/\kc_M$. Classically,
(\ref{ases}) follows from the equality between the norm, modulo compact operators, of a singular integral operator on a compact manifold
and the supremum norm of its symbol, proven by Gohberg \cite{Goh}\ and Seeley \cite{See}.
Proofs of that result in the language of pseudodifferential operators appeared in \cite{Hp,KN}.


\section{Preliminaries}

In this section, we define the C*-algebra $\ac$ and review
results about its structure that will be needed in the next and final section.

Let $B$ denote a compact Riemannian manifold of dimension $n$ without boundary. Consider the cylinder
$\Omega = \mathbb{R}\times B$ and denote by $\Delta_{\Omega}=\Delta_\R+\Delta_B$ its Laplacian 
(we use the classical analysis convention sign for the Laplacian, $-\Delta_{\Omega}\geq 0$)
Define $\ac$ as the C*-subalgebra of ${\bc}(L^2(\Omega))$
 (we denote by $\bc({\mathcal H})$ the C*-algebra of
all bounded operators on a Hilbert space ${\mathcal H}$) generated by: 
\begin{itemize}
\item[(i)] all $A_1=a(M_x)$, operators of multiplication by $ a\in C^\infty({B})$,
  \[
    [a(M_x)f](t,x) = a(x)f(t,x),\ \ (t,x)\in \mathbb{R}\times B;
    \]
\item[(ii)] all $A_2=b(M_t)$, operators of multiplication by $b\in C([-\infty,+\infty])$ (this denotes the set of all
  continuous functions on $\mathbb{R}$ which have limits at $-\infty$ and $+\infty$),
  \[
    [b(M_t)f](t,x) = b(t)f(t,x),\ \ (t,x)\in \mathbb{R}\times B;
    \]
  \item[(iii)] all operators $A_3=e^{ikM_{t}}$, $k\in\Z$, defined by
    \[
    [e^{ikM_{t}}f](t,x) = e^{ikt}f(t,x),\ \ (t,x)\in \mathbb{R}\times B;
    \]  \label{generators}  
\item[(iv)] $A_4=\Lambda:=(1-\Delta_{\Omega})^{-1/2};$ 
\item[(v)] $A_5=\frac{1}{i}\frac{\partial}{\partial t}\Lambda;$
\item[(vi)] $D_x\Lambda$, where $D_x=-i\sum_{k=1}^{n}c^k(x)\partial_{x_k}$ is a smooth vector field on $B$. 
\end{itemize}

Since $\Omega$ is a complete manifold, its Laplacian, defined on $C_c^\infty(\Omega)$, has a unique extension
(also denoted by $\Delta_\Omega$) as an unbounded self-adjoint operator on $L^2(\Omega)$
(\cite{esa}, \cite[Theorem IV.1.8]{CO1}) and one can then define $\Lambda$ by functional calculus
. 
The operators in (v) and (vi), a priori defined on the dense subset $\Lambda^{-1}(C_c^\infty(\Omega))$,
have continuous extensions to $L^2(\Omega)$ \cite[Section V.1]{CO1}.

Moreover, if $L$ denotes an $M$-th order differential operator on $\Omega$
whose (smooth) coefficients approach $2\pi$-periodic functions as $t$ tends to $-\infty$ and
as $t$ tends to $+\infty$,  the operator $L\Lambda^M$
can be extended to a bounded operator on $L^2(\Omega)$, which belongs to $\ac$ . These $L\Lambda^M$ are 
zero-order {\em classical} pseudodifferential operators (meaning that their local symbols have asymptotic expansions in homogeneous components). 

Applying Cordes' comparison algebra techniques and using results from \cite{CO3,CM},
the statements about the structure of $\ac$ contained in the rest of this Section were proven in \cite{T2} (see also \cite{MTh}).

Let $\ec$ denote the commutator ideal of $\ac$, i.e, the smallest closed ideal containing all commutators $AB-BA$, $A$ and $B$ in $\ac$. 
Let $\overline{\Omega}$ denote the compactfication $[-\infty,+\infty]\times B$ of $\Omega$, and let $S^*\!M$ denote the cosphere bundle of 
any manifold $M$. We may then write 
\[
S^*\overline{\Omega}=\{(t,x,\tau,\xi);\,(t,\tau)\in T^*[-\infty,+\infty],\,(x,\xi)\in T^*B,\,|\tau|^2+|\xi|^2=1\}.
\]

\begin{Teo}\label{ma2} The spectrum of the commutative C*-algebra $\ac/\ec$ is homeomorphic to the following subset of 
$(S^*\overline{\Omega})\times S^1:$
\[
\mathbf{M}=\{(t,x,\tau ,\xi,e^{i\theta});\ \theta\in\R\ \mbox{and}\ \theta=t\ \mbox{if}\ |t|<\infty\}. 
\]
Denoting by
\begin{equation}
  \label{sigmasymb}
\sigma:\ac\to C(\mathbf{M})
\end{equation}
the composition of the Gelfand isomorphism $\ac/\ec\simeq C(\mathbf{M})$ with the projection $\ac\to\ac/\ec$,
for each of the generators of $\ac$ listed above,
$\sigma(A_k)(t,x,\tau ,\xi,e^{i\theta})$, $k=1,\cdots,6$, is equal to, respectively, 
\begin{equation}\label{sigma}
 a(x),~ ~ ~ b(t),~ ~ ~ e^{ij\theta},~ ~ ~ 0, ~ ~ ~ \tau,~ ~ ~ \sum_{k=1}^nc^j(x)\xi_j.
\end{equation}
\end{Teo}
Since $M$ is a compactification of $S^*\Omega$, we may regard $C(\mathbf{M})$ as a subspace of $C(S^*\Omega)$. Theorem~\ref{ma2}\ can then be rephrased by
saying that the usual principal symbol of classical pseudodifferential operators, defined on a dense subalgebra of $\ac$, 
has a continuous extension $\sigma$ to $\ac$ whose image is equal to $C(\mathbf{M})$ and whose kernel is $\ec$.

Let $F_d: L^2(S^1)\rightarrow {\ell}^2(\mathbb{Z})$ denote the discrete Fourier transform
$$
(F_d u)_j = \int_{0}^{1} u(e^{2\pi i\varphi}) e^{-2\pi ij\varphi} d\varphi ~,~ u\in L^2(\mathbb{S}^1) \mbox{ and	 } j\in \mathbb{Z} 
$$
($S^1$ is equipped with the Lebesgue measure pushforwarded by $[0,1)\ni\varphi\mapsto e^{2\pi i\varphi}\in S^1$; hence $F_d$ is unitary).
For each sequence $(b_j)_{j\in\Z}$ with limits as $j\to -\infty$ and as $j\to +\infty$, we define bounded operators $b(M_j)$ on $\ell^2(\Z)$ 
and $b(D_\varphi)$
on $L^2(S^1)$ by
\[
[b(M_j)u]((u_j)_{j\in\Z})=(b_ju_j)_{j\in\Z}\ \ \ \mbox{and}\ \ \ b(D_\varphi)=F_d^{-1}b(M_j)F_d.
\]
Given $a\in C(S^1)$, let $a(M_\varphi)\in\bc(L^2(S^1))$ denote the operator of multiplication by $a$. Finally let $\Sc$
denote the C*-subalgebra of ${\mathcal{L}}(L^2({S}^1))$ generated by the operators of the form $a(M_\varphi)$ and $b(D_\varphi)$, for $a\in C(S^1)$
and $b=(b_j)_{j\in\Z}$ a sequence as above.

In Cordes' language, $\Sc$ is the only comparison algebra over $S^1$. The principal symbol,
defined on the dense *-subalgebra of all classical zero-order pseudos, extends to a surjective C*-algebra homomorphism
(see the Introduction)
\begin{equation}
\label{sigmacircle}
\sigma_S:\Sc\to C(S^1\times\{-1,+1\}),
\end{equation}
whose kernel is equal to the set  $\kc_{S^1}$ of all compact operators on $L^2(S^1)$. On our set of generators,
one has $[\sigma(a(M_\varphi))](e^{2\pi i\varphi},\pm 1)=a(e^{2\pi i\varphi})$, for $a\in C(S^1)$, and $[\sigma(b(D_\varphi))](e^{2\pi i\varphi},\pm 1)=b(\pm\infty)$,
  for $b=(b_j)_{j\in\Z}$ a sequence having limits as $j\to +\infty$ and as $j\to -\infty$. A detailed exposition of some of these facts can be
  found in \cite{R}.

  For each  $\varphi \in \mathbb{R}$, let $U_{\varphi}$ be the operator on $L^2(S^1)$ given by
  $U_{\varphi}f (z) = z^{-\varphi}f(z)$, $z\in S^1$ (we take
  the principal branch of $z^{-\varphi}=e^{-\varphi \log z}$, $z\neq 0$).
  Let $Y_\varphi := F_d U_{\varphi} F_d^{-1}$. 
This defines a smooth family  $\varphi\mapsto Y_{\varphi}$ of unitary operators on $\ell^2(\mathbb{Z})$ such that for all 
$k\in \mathbb{Z},~(Y_k u)_j = u_{j+k},$ and $Y_{\varphi} Y_{\omega} = Y_{\varphi +\omega}$.
Given $u\in L^2(\R)$, the sequence $\tilde{u}(\varphi)=(u(\varphi-j)_{j\in\Z})$ belongs to
  $\ell^2(\Z)$ for almost every $\varphi\in[0,1)$ and the map
  \[
  \begin{array}{rcl}
    W:L^2(\R)&\longrightarrow &L^2(S^1)\otimes\ell^2(\Z)\simeq L^2(S^1;\,\ell^2(\Z))\\
    u&\longmapsto&Wu,\ \ \ \ \ \, (Wu)(e^{2\pi i\varphi})=Y_\varphi\tilde{u}(\varphi)
\end{array}
    \]
    is a Hilbert space isomorphism.

    We have used $\otimes$ to denote Hilbert-space tensor product. Given $A\in\mathfrak{B}(\mathcal{H})$ and
    $B\in\mathfrak{B}(\mathcal{K})$, we also denote by $A\otimes B$ the bounded operator on $\mathcal{H}\otimes\mathcal{K}$ defined
    by $(A\otimes B)(u\otimes v)=Au\otimes Bv$. And also, given C*-algebras $\mathfrak{A}_H\subset\mathfrak{B}(\mathcal{H})$ and
$\mathfrak{A}_K\subset\mathfrak{B}(\mathcal{K})$, we
      denote by $\mathfrak{A}_H\otimes\mathfrak{A}_K\subset\mathfrak{B}(\mathcal{H}\otimes\mathcal{K})$ the closed linear span
      of all $A\otimes B$, $A\in\mathfrak{A}_H$ and $B\in\mathfrak{A}_K$.

We will denote by  $\kc_X$ denote the sets of all compact
operators on $L^2(X)$, for $X=\Omega$, $S^1$, $\Z$, $B$, $\Z\times B$, $S^1\times B$.
We will also write $W$ and $F$ meaning $W\otimes I_B$ and
$F\otimes I_B$, where $I_B$ denotes the identity map on $L^2(B)$. Similarly we will also denote by $\sigma_S$ the
map $\sigma_S\otimes I_{\kc_{\Z}}\otimes I_{\kc_{B}}$ defined on $\Sc\otimes\kc_\Z\otimes\kc_B$.
  
      The *-isomorphism
        $\mathfrak{B}(L^2(\Omega))\ni A\ \longmapsto\ WF^{-1}AFW^{-1}\in\mathfrak{B}(L^2(S^1)\otimes\ell^2(\Z)\otimes L^2(B))$
        restricted to $\ec$ defines a C*-isomorphism
\begin{equation}\label{e}
                      \ec\       \xrightarrow{\cong}   \ \Sc\otimes\kc_\Z\otimes\kc_B.
\end{equation}
  The kernel of the surjective map
        $\sigma_S:\Sc\otimes\kc_\Z\otimes\kc_B\to C(S^1\times\{-1,+1\})\otimes\kc_\Z\otimes\kc_B$ is equal to
        $\kc_{S^1}\otimes\kc_\Z\otimes\kc_B\cong\kc_\Omega$ and hence $\sigma_S$ induces a C*-algebra isomorphism 
  \begin{equation}\label{psi}
        \ec/\kc_\Omega \xrightarrow{\cong} C(S^1\times\{-1,+1\})\otimes\kc_\Z\otimes\kc_B\ \cong\
        C(S^1\times\{-1,+1\};\,\kc_\Z\otimes\kc_B).
        \end{equation}
  The composition of this isomorphism with the canonical projection $\ec\to\ec/\kc_\Omega$ defines a C*-algebra homomorphism
        \begin{equation}\label{gsymbol}
\gamma:\,\ec\ \longrightarrow\ C(S^1\times\{-1,+1\};\,\kc(L^2(\Z\times B))),
        \end{equation}
        which can be extended, via the left regular representation, to a map from $\ac$ 
        to the algebra of bounded operators on the Banach space $C(S^1\times\{-1,+1\};\,\kc(L^2(\Z\times B))$.
        The image of this Banach-algebra homomorphism is contained
        in $C(S^1\times\{-1,+1\};\,\Bc(L^2(\Z\times B))$ (an $f\in C(S^1\times\{-1,+1\};\bc(L^2(\Z\times B))$ is here regarded
        as a multiplication operator on $C(S^1\times\{-1,+1\};\kc(L^2(\Z\times B))$\,) and this defines a C*-algebra homomorphism, which we also
        denote by $\gamma$,
        \begin{equation}
          \label{gamma}
          \gamma:\,\ac\ \longrightarrow\ C(S^1\times\{-1,+1\};\,\bc(L^2(\Z\times B))).
        \end{equation}
We define yet another *-homomorphism 
\begin{equation}
          \label{Gamma}
\Gamma:\,\ac\ \longrightarrow\ C(S^1\times\{-1,+1\};\,\bc(L^2(S^1\times B)))
\end{equation}
\[
\Gamma_A(e^{2\pi i\varphi},\pm 1)=F_d^{-1}\Gamma_A(e^{2\pi i\varphi},\pm 1)F_d,\ \ e^{2\pi i\varphi}\in S^1.
\]
The image of $\Gamma$ is contained in $C(S^1\times\{-1,+1\};\,\Psi(S^1\times B))$, with $\Psi(S^1\times B)$ as defined at
the end of the Introduction.  
It follows almost by definition that the image of $\ec$ under $\Gamma$ equals
$C(S^1\times\{-1,+1\};\,\kc_{S^1\times B})$.

In order to give a more explicit description of $\Gamma$, we need two more definitions.
For each $\varphi\in\R$, let us denote by $e^{i\varphi M_\theta}$ the multiplication operator (by a possibly discontinuous function)
\[
[e^{i\varphi M_\theta}u](e^{i\theta})=e^{i\varphi \theta}u(e^{i\theta} ),\ \ \ -\pi<\theta\leq\pi,\ \ \ u\in L^2(S^1).
\]
And let us denote by $D_\theta$ the densely defined operator on $L^2(S^1)$
\[
(D_\theta u)(e^{i\theta})=-i\tilde{u}'(\theta),\ \ \ \tilde{u}(\theta)=u(e^{i\theta}),\ \ \       \theta\in\R,    \ \ \ u\in C^\infty(S^1).
\]
For each of the generators of $\ac$ defined on page~\pageref{generators}, the value of $\Gamma_{A_k}(e^{2\pi i\varphi},\pm 1)$, $1\leq k\leq 6$, equals
\begin{equation}\label{explicito}
a(M_x),\ \ \ b(\pm\infty)I,\ \ \ e^{ijM_\theta},\ \ \ e^{-i\varphi M_\theta}[1+(D_\theta-\varphi)^2-\Delta_B]^{-1/2}e^{i\varphi M_\theta},
\end{equation}
\[ e^{-i\varphi M_\theta}(D_\theta-\varphi)[1+(D_\theta-\varphi)^2-\Delta_B]^{-1/2} e^{i\varphi M_\theta},
\]\[
e^{-i\varphi M_\theta}D_x[1+(D_\theta-\varphi)^2-\Delta_B]^{-1/2} e^{i\varphi M_\theta},
\]
respectively (the first three classes of generators are mapped to constant functions).

For any positive integer $k$, the homomorphisms $\sigma$ and $\Gamma$
defined in (\ref{sigmasymb}) and (\ref{Gamma}) canonically induce homomorphisms, which we also denote 
by $\sigma$ and $\Gamma$, on the algebra $M_k(\ac)$ of $k$-by-$k$ matrices with entries in $\ac$. Regarding
$M_k(\ac)$ as a C*-subalgebra of $\mathfrak{B}(L^2(\Omega)^k)$, it turns out that an $A\in M_k(\ac)$ is a Fredholm operator
if and only if $\sigma(A)$ and $\Gamma(A)$ are invertible \cite[Theorem 3.2]{T2}. Moreover, the quotient of $M_k(\ac)$ by the ideal of
compact operators on $L^2(\Omega)^k$ is isomorphic to the image of the {\em total symbol}
\begin{equation}
  \label{totalsymbol}
\sigma\oplus\Gamma: M_k(\ac)\ \longrightarrow\ M_k(C(\mathbf{M}))\oplus M_k(C(S^1\times\{-1,+1\};\,\Psi(S^1\times B))).
\end{equation}


\section{The principal-symbol exact sequence}

In this section we analyse the $K$-theory index map associated to the exact sequence of C*-algebras
\begin{equation}
  \label{1}
0\ \longrightarrow\ \ec\ \hookrightarrow\ \ac\ \stackrel{\sigma}{\longrightarrow}\ C(\mathbf{M})\ \longrightarrow\ 0,
\end{equation}
induced by the continuous extension to $\ac$ of the usual principal symbol of pseudodifferential operators (see Theorem~\ref{ma2}).

We first consider the exact sequence of C*-algebras associated to the map $\sigma_S$ defined in (\ref{sigmacircle}),
\begin{equation}
\label{circleexact}
0 \ \longrightarrow \kc_{S^1} \  \hookrightarrow \ \Sc \  \stackrel{\sigma_S}{\longrightarrow}
\ C(S^1\times\{-1,+1\}) \ \longrightarrow \ 0.
\end{equation} 
This sequence tells us that an $A\in\Sc$ is a Fredholm operator if and only if $\sigma_S(A)$ never vanishes. A well-known 
index formula for zero-order pseudodifferential operators on the circle, which essentially goes back to 
Noether \cite{No} (see also the example right after Theorem~19.2.4 in \cite{H3}), implies that the index of a 
Fredholm operator $A\in\Sc$ is equal to the difference of the winding numbers of the restrictions of $\sigma_S(A)$ 
to the two copies of $S^1$ in $S^1\times\{-1,+1\}$. A detailed proof of this fact without explicit mention of 
pseudodifferential or singular integral operators can be found in \cite[Se\c{c}\~ao 10]{R}. It then follows that 
there is a Fredholm operator of index 1 in $\Sc$ and, therefore, the index map $\delta_1$ in the 
standard exact sequence of $K$-groups associated to (\ref{circleexact})
\[
\begin{array}{ccccc}
 \mathbb{Z} \cong K_0(\kc_{S^1})  & \! {\longrightarrow} & \! 
 K_0(\Sc)  & \! {\longrightarrow} & \!
 K_0(C(S^1\times\{-1,+1\}))  \\ \\    
 {\delta}_1 \ \uparrow & \! ~  & \! ~  & \! ~ & \! \downarrow \ {\delta}_0  \\ \\
 K_1(C(S^1\times\{-1,+1\})) & \! {\longleftarrow} & \! 
 K_1(\Sc)  & \! {\longleftarrow} & \!
 K_1(\kc_{S^1}) = 0
\end{array}
\]
is surjective. We then get $K_0(\Sc)\cong\Z^2$ and $K_1(\Sc)\cong\Z$. 

\begin{Pro}\label{prop1}
The canonical projection $\pi:\ec\to\ec/\kc_\Omega$ induces a $K_0$-isomorphism. 
\end{Pro} 
\Dem The previous considerations and the isomorphisms (\ref{e}) and (\ref{psi}) imply that
\[
K_0(\ec)\cong K_0(\ec/\kc_\Omega)\cong K_1(\ac/\ec)\cong\Z^2 \ \ \ \mbox{and}\ \ \ K_1(\ec)\cong\Z.
\]
The standard exact sequence of $K$-groups associated to 
\[
0 \ \longrightarrow \kc_{\Omega} \  \hookrightarrow \ {\ec} \  \stackrel{\pi}{\longrightarrow}
\ \frac{\ec}{\kc_{\Omega}} \ \longrightarrow \ 0
\]
then becomes
\begin{equation}\label{six}
\begin{array}{ccccc}
 \mathbb{Z}& \! {\longrightarrow} & \! 
 \Z^2 & \! {\longrightarrow} & \!
 \Z^2  \\ \\    
 {\delta}_1 \ \uparrow & \! ~  & \! ~  & \! ~ & \! \downarrow \ {\delta}_0  \\ \\
 \Z^2 & \! {\longleftarrow} & \! 
 \Z  & \! {\longleftarrow} & \!
 0
\end{array}
\end{equation}
The index map $\delta_1$ is non-zero (otherwise there would be an exact sequence $0\to\Z\to\Z^2\to\Z^2\to 0$).
If it were not surjective, the image of the upper left homomorphism
would be a finite non-trivial group, which could not be a subgroup of $\Z^2$. The exactness of the
sequence then provides the desired result.
\cqd

The following rephrasing of the fact that $\delta_1$ in (\ref{six}) is surjective is perhaps interesting in itself.
\begin{Cor}
  There exists $k>0$ and a Fredholm operator of index $1$ in $M_k(\ec)\subset\Bc(L^2(\Omega)^k)$.\end{Cor}

We will need to use Banach algebra $K$-theory \cite{BB}, in which $K_0$ consists of formal differences of classes of idempotents,
while the elements of $K_1$ are classes of invertibles. The following lemma was used before in \cite{Hess,TO1,N}, 
see \cite[Lemma 1]{TO1}\ for a proof.

\begin{Lem}\label{lema1} Let $\mathcal{A}$ be a unital Banach algebra, $\mathcal{J}$ a closed ideal of $\mathcal{A}$,
  and let $\delta_1: K_1(\mathcal{A}/\mathcal{J})\rightarrow K_0(\mathcal{J})$ denote the index map
  associated to the exact sequence 
  $0\rightarrow \mathcal{J} \rightarrow \mathcal{A} \rightarrow \mathcal{A}/\mathcal{J} \rightarrow 0$.
  If $u\in M_k(\mathcal{A}/\mathcal{J})$ is invertible, $\pi(a)=u$ and $\pi(b)=u^{-1},$ then 
\[
\delta_1([u]_1)= 
\left[
\left(\begin{array}{cc}
2ab-(ab)^2 & a(2-ba)(1-ba) \\
(1-ba)b & (1-ba)^2 \\
\end{array}\right)
\right]_0
-
\left[
\left(\begin{array}{cc}
1 & 0 \\
0 & 0 \\
\end{array}\right)
\right]_0
\]
\end{Lem}

Next we consider the exact sequence obtained from (\ref{1}) quotienting $\ec$ and $\ac$ by the compact ideal,
\begin{equation}\label{2}
0 \ \longrightarrow \ \frac{\ec}{{\kc}_{\Omega}} \  
\longrightarrow \ \frac{\ac}{{\kc}_{\Omega}} \  
\stackrel{\pi}{\longrightarrow} \ \frac{\ac}{\ec} \ \longrightarrow \ 0.
\end{equation}

The following Proposition describes the connecting map $\delta_1:K_1(\ac/\ec)\to K_0(\ec/\kc_\Omega)$ in terms of
the index of Fredholm operators on $L^2(S^1\times B)$. We recall that each element of $K_1(\ac/\ec)$ is of the form
$[[A]_{\ec}]_1$, where $A\in M_k(\ac)$ is invertible modulo $M_k(\ec)$
for some positive integer $k$ and $[A]_{\ec}$ denotes its class in the quotient $M_k(\ac)/M_k(\ec)$. 

\begin{Pro}\label{delta1} Let $A \in M_k(\ac)$ be  such that $[A]_{\ec}$ is invertible in  $M_k(\ac/\ec)$ for some positive integer $k$.
  Then  ${\Gamma}_{A}(z,-1)$ and  ${\Gamma}_{A}(z,+1)$ are Fredholm operators on $L^2(S^1\times B)^k$ for every $z \in S^1$ and
\begin{equation}\label{ph}
{\delta}_{1}(\left[[A]_{\ec}]_1\right) = 
( \mathsf{ind}\,{\Gamma}_{A}\!(1,-1)\,[E]_0,  \mathsf{ind}\,{\Gamma}_{A}\!(1,+1)\,[E]_0)
\end{equation}
where $E$ is a rank one projection on $L^2(S^1\times B)$     
and $\mathsf{ind}$ denotes the Fredholm index $($\!It will become clear in the proof how the left side of $(\ref{ph})$ can
be regarded as an element of $K_1(\ec/\kc_\Omega))$. 
\end{Pro}

\Dem Recall that the image of $\ec$ under $\Gamma$ equals
$C(S^1\times\{-1,+1\};\,\kc_{S^1\times B})$.

Let $B\in M_k(\ac)$ be such that $I-AB$ and $I-BA$ belong to $M_k(\ec)$. Then, $I-\Gamma_A\Gamma_B$ and $I-\Gamma_B\Gamma_A$ belong to
$C(S^1\times\{-1,+1\};\,M_k(\kc_{S^1\times B}))$. This means that, at each $(z,\pm 1)$ in $S^1\times\{-1,+1\}$,
$\Gamma_A(z,\pm 1)$ is invertible modulo compacts, and hence it is a Fredholm operator. 

Let $a=[A]_{\kc}$ and $b=[B]_{\kc}$. By Lemma~\ref{lema1}, 
\begin{equation}\label{vn}
\delta_1([[A]_{\ec}]_1)= 
\left[
\left(\begin{array}{cc}
2ab-(ab)^2 & a(2-ba)(1-ba) \\
(1-ba)b & (1-ba)^2 \\
\end{array}\right)
\right]_0
-
\left[
\left(\begin{array}{cc}
1 & 0 \\
0 & 0 \\
\end{array}\right)
\right]_0
\in K_0(\ec/\kc_\Omega)
\end{equation}
(we are regarding $\ec/\kc_\Omega$ as an ideal of $\ac/\kc_\Omega$). Using the isomorphism (\ref{psi}) and the definitions of $\gamma$ and $\Gamma$
that follow (\ref{psi}), we see that the map
\[
\ec/\kc_\Omega\ni [T]\ \ \longmapsto\ \ \Gamma_T\in C(S^1\times\{-1,+1\};\,\kc_{S^1\times B})
\]
is an isomorphism. Using this isomorphism as an identification, we may replace $a$ and $b$ by $\Gamma_A$ and $\Gamma_B$ in (\ref{vn}) and get
$
\delta_1([[A]_{\ec}]_1)= 
$
\[
\left[
\left(\begin{array}{cc}
2\Gamma_A\Gamma_B-(\Gamma_A\Gamma_B)^2& \Gamma_A(2-\Gamma_B\Gamma_A)(1-\Gamma_B\Gamma_A) \\
(1-\Gamma_B\Gamma_A)\Gamma_B & (1-\Gamma_B\Gamma_A)^2 \\
\end{array}\right)
\right]_0
-
\left[
\left(\begin{array}{cc}
1 & 0 \\
0 & 0 \\
\end{array}\right)
\right]_0
\]\[
\in K_0(C(S^1\times\{-1,+1\};\,\kc_{S^1\times B}))\cong
K_0(C(S^1;\,\kc_{S^1\times B}))  \oplus
K_0(C(S^1;\,\kc_{S^1\times B}))
\]
Let us consider now the map $p:C(S^1;\,\kc_{S^1\times B})\to\kc_{S^1\times B}$, $p(f)=f(1)$. The map $q:\kc_{S^1\times B}\to C(S^1;\,\kc_{S^1\times B})$
that regards a compact operator as a constant function is a right inverse for $p$. Moreover, the kernel of $p$ is the suspension of $\kc$, and then
its $K_0$-group is 0. All that imply that $p$ induces a $K_0$-isomorphism. Using this isomorphism as an identification, we may write
$\delta_1([[A]_{\ec}]_1)$ as the element of $K_0(\kc_{S^1\times B})\oplus K_0(\kc_{S^1\times B})$
\[
\left(
\left[
\left(\begin{array}{cc}
2A^-B^--(A^-B^-)^2& A^-(2-B^-A^-)(1-B^-A^-) \\
(1-B^-A^-)B^- & (1-B^-A^-)^2 \\
\end{array}\right)
\right]_0
-
\left[
\left(\begin{array}{cc}
1 & 0 \\
0 & 0 \\
\end{array}\right)
\right]_0
\right.,
\]
\[
\left.
\left[
  \left(\begin{array}{cc}
2A^+B^+-(A^+B^+)^2& A^+(2-B^+A^+)(1-B^+A^+) \\
(1-B^+A^+)B^+ & (1-B^+A^+)^2 \\
   \end{array}\right)
\right]_0
-
\left[
\left(\begin{array}{cc}
1 & 0 \\
0 & 0 \\
\end{array}
\right)
\right]_0
\right),
\]
where we have denoted $\Gamma_A(1,\pm 1)$ by $A^\pm$ and analogously for $\Gamma_B(1,\pm 1)$.

If we apply Lemma \ref{lema1} to $\mathcal{A}=\mathcal{B}(L^2(S^1\times B))$, $\mathcal{J}=\kc_{S^1\times B}$ and $u=[A^\pm]_{\kc}$, recalling that
  the K-theory index map is a generalization of the Fredholm index \cite{RLL}, we conclude that this element of
  $K_0(\kc_{S^1\times B})\oplus K_0(\kc_{S^1\times B})$ that we have just written above is equal to the left side of (\ref{ph}). \cqd

  Theorem~\ref{thm2} below follows immediately from Proposition~\ref{prop1} and Proposition~\ref{delta1}.
  To simplify notation, in its statment we regard as identifications  the canonical group isomorphisms
  \[
  K_0(C(S^1\times\{-1,+1\};\kc_{S^1\times B}))\cong
  K_0(C(S^1;\kc_{S^1\times B}))\oplus K_0(C(S^1;\kc_{S^1\times B}))\cong \Z\oplus\Z.
  \]

  \begin{Teo}\label{thm2}
    Let $\delta_1^/:K_1(\ac/\kc)\to K_0(\ec)$ denote the index map associated to the exact sequence
    $0\to\ec\to\ac\to\ac/\ec\to 0$, let $\Gamma_*:K_0(\ec)\to\Z\oplus\Z$ denote the isomorphism induced by the surjective map
$\Gamma:\ec\to C(S^1\times\{-1,+1\};\kc_{S^1\times B})$,
and denote $\delta_1=\Gamma_*\circ\delta_1^/$. Given $A \in M_k(\ac)$  such that $[A]_{\ec}$ is invertible in  $M_k(\ac/\ec)$ for some positive integer $k$,
then  $A^-={\Gamma}_{A}(z,-1)$ and  $A^+={\Gamma}_{A}(z,+1)$ are Fredholm operators on $L^2(S^1\times B)^k$ and
\begin{equation}\label{quase}
  {\delta}_{1}(\left[[A]_{\ec}]_1\right) = 
( \mathsf{ind}A^-,\mathsf{ind}A^+)
\end{equation}
where $\mathsf{ind}$ denotes the Fredholm index.
  \end{Teo}
  
  The final goal should be to describe the index map associated to $(\ref{1})$ in purely topological terms.
  That can now be easily achieved if we invoke the Atiyah-Singer Theorem.

  In the following we adopt topological K-theory definitions and notation of \cite[Chapter I]{BB}. We denote, as usual, by $T^*M$, $B^*M$ and $S^*M$,
  respectively, the tangent bundle, the bundle of unit closed balls in the tangent bundle, and the co-sphere bundle of any manifold $M$.
 
  \begin{Teo}\label{topology} Let $\delta_1$ denote the index map associated to $(\ref{1})$, and let $f\in M_k(C(\mathbf{M}))$ be invertible for some
    positive $k$. The restrictions $f_-$ and $f_+$ of $f$ to the points of $\mathbf{M}$ with, respectively,  $t=-\infty$ and $t=+\infty$ are invertible elements of
    $M_k(C(S^*(S^1\times B)))$. Let
    \[
    x_\pm=[(E^k,E^k,f_\pm)]\in K(B^*(S^1\times B),S^*(S^1\times B))\cong K(T^*(S^1\times B))
    \]
    denote the classes of the triples $(E^k,E^k,f_\pm)$,
    where $E^k$ denotes the trivial bundle
    $(S^1\times B)\times \C^k$. Then
    \[
    \delta^1([f]_1)=(\mathsf{ind}_t(x_+),\mathsf{ind}_t(x_-))\in\Z\oplus\Z,
    \]
    where $ind_t:K(T^*(S^1\times B))\to\Z$ denotes the topological index of Atiyah and Singer \cite{AS}.
    \end{Teo}

  \Dem
  Let $A\in M_k(\ac)$ be such that $\sigma(A)=f$. Define $A^\pm=\Gamma_{A}(1,\pm 1)$. We claim that $A^-$ and $A^+$ belong to $M_k(\Psi(S^1\times B))$ and
  $\sigma_{S^1\times B}(A^\pm)=f_\pm$. By definition, $f_\pm=\sigma(A)|_{t=\pm\infty}$. The maps $\ac\ni A\mapsto \sigma(A)|_{t=\pm\infty}$,
  $\ac\ni A\mapsto \Gamma_{A}(1,\pm 1))|_{t=\pm\infty}$ and
  $\Psi(M)\ni T\mapsto \sigma_{S^1\times B}(T)$ are C*-algebra homomorphisms. So, it is enough to verify our claim on a set of generators of $\ac$. 
  For each of the generators $A_1,A_2,\cdots,A_6$  defined on page~\pageref{generators}, the claim follows immediately from (\ref{explicito}).

  Our theorem would already follow from the Atiyah-Singer Theorem \cite{AS} if $A^+$ and $A^-$ were classical pseudodiferential operators. For the general case,
  we need to show that the topological index of the class $[(E^k,E^k,\sigma(T))]\in K(T^*(S^1\times B))$ is equal to the Fredholm index of any Fredholm
  operator $T\in\Psi(S^1\times B)$ (the Atiyah-Singer Theorem gives that for those $T$ that belong to the dense subalgebra of all classical
  pseudodifferential operators). 

In the rest of this proof, we will replace $S^1\times B$ by an arbitrary compact manifold $M$.

  Let $\delta^{tbs}$ denote index mapping for the exact sequence of C*algebras
  \[
0\longrightarrow C_0(T^*M)\longrightarrow C(B^*M)\longrightarrow C(S^*M)\longrightarrow 0
  \]
  and let $\iota:K_0(C_0(T^*M))\mapsto K(T^*M)$ be the canonical isomorphism between these two groups.
  We then have
  \begin{equation}
    \label{tbs}
\iota\circ\delta^{tbs}([g]_1)=[(E^k,E^k,g)], \ \ \mbox{for any invertible}\ \ g\in M_k(C(S^*M)).
\end{equation}
For a proof of (\ref{tbs}) in a slightly more general setting, see \cite[Proposition 15]{MNS}, for example.

  Now let $T\in M_k(\Psi(M))$ be a Fredholm operator. Its Fredholm index $\mathsf{ind}(T)$ coincides with $\delta^M([\sigma_M(T)]_1)$,
  where $\delta^M$ denotes the index map of (\ref{ases}). Because every $f\in C^\infty(S^*M)$ is the symbol of a zero-order classical
  pseudodifferential operator, there exists such a pseudo $T_0$ with invertible symbol so close to $T$ that $[\sigma(T)]_1=[\sigma(T_0)]_1$. Then
  \[
  \mathsf{ind}(T)=\delta^M([\sigma_M(T)]_1)=\delta^M([\sigma_M(T_0)]_1)=\mathsf{ind}_t([(E^k,E^k,\sigma_M(T_0))])=
  \]
  \[
\mathsf{ind}_t(\iota\circ\delta^{tbs}([\sigma_M(T_0)]_1))=\mathsf{ind}_t(\iota\circ\delta^{tbs}([\sigma_M(T)]_1))=\mathsf{ind}_t([(E^k,E^k,\sigma_M(T_0))]).
  \]
 The Atiyah-Singer Theorem was used in the third of the above equalities. Twice we used (\ref{tbs}). \cqd
  
 {\bf Example}

Some statements made here about operators in $\ac$ are proven in \cite{T2}, using results from \cite{CO1}.
 
 Let $L:C^\infty(\Omega)^k\to C^\infty(\Omega)^k$ be a differential operator of the form
 \[
 L=\sum_{j+|\alpha|\leq N}a_{j,\alpha}(t,x)
 \left(\frac{1}{i}\frac{\partial}{\partial t}\right)^j
 \left(\frac{1}{i}\frac{\partial}{\partial x}\right)^\alpha,
 \]
 where each $a_{j,\alpha}\in M_k(C^\infty(\Omega))$ has support contained in $U\times\R$, $U$ the domain of a chart in $B$. Moreover, let us assume that
 $a_{j,\alpha}$ is {\em semi-periodic} in the following sense. There exist $2\pi$-periodic (in $t$) $a_{j,\alpha}^\pm\in M_k(C^\infty(\Omega))$ and,
 for a choice of real-valued $\chi^\pm\in C^\infty(\R)$ such that $\chi^+(t)+\chi^-(t)=1$ for all $t$ and $\chi_\pm(t)=0$ for $\mp t>1$, there
 exists $a^0_{j,\alpha}\in M_k(C^\infty(\Omega))$ vanishing at infinity such that
 \[
a_{j,\alpha}(t,x) = \chi^+(t)a_{j,\alpha}^+(t,x) + \chi^-(t)a_{j,\alpha}^-(t,x) + a^0_{j,\alpha}(t,x), \ \ \mbox{for all}\ \ (t,x)\in\Omega.
\]

The operator $A=L\Lambda^N$\!, a priori defined on the dense subspace $[\Lambda^{-N}(C_c^\infty(\Omega))]^k$ of $L^2(\Omega)^k$,
extends to a bounded operator $A\in M_k(\ac)$ (\cite[Sections VII-3 and IX-3]{CO1}, \cite[Theorem 3.7]{T2}).
Each of the generators $A_1,\cdots,A_5$ defined on page~\pageref{generators}
is a particular example of such an $A$, with $k=1$ and $N$ equal to 0 or 1. And because $B$ has a finite atlas, each $A_6$ is
a finite sum of such $A$'s.

The symbol $\sigma_A$ coincides with the continuous extension to the compactification
$\mathbf{M}$ of $S^*\Omega$ of the principal symbol of $L$,
\[
\sigma(A)(t,x,\tau,\xi)=\sum_{j+|\alpha|=N}a_{j,\alpha}(t,x)\tau^j\xi^\alpha,\ \  (t,x,\tau,\xi)\in S^*\Omega
\]
It follows that $A$ is invertible modulo $\ec$ if and only if $L$ is {\em uniformily elliptic} in the sense that
\[
\inf\{\,|\!\!\!\sum_{j+|\alpha|= N}a_{j,\alpha}(t,x)\tau^j\xi^\alpha|,\ (t,x,\tau,\xi)\in S^*\Omega\}\ >\ 0.
\]
We also have 
\[
\Gamma_A(1,\pm 1)= \sum_{j+|\alpha|=N}a_{j,\alpha}^\pm(t,\theta)
 \left(\frac{1}{i}\frac{\partial}{\partial \theta}\right)^j
 \left(\frac{1}{i}\frac{\partial}{\partial x}\right)^\alpha\Lambda^N
 \]
 ($\theta$ denotes the argument of a point in $S^1$). These are zero-order elliptic classical pseudodifferential
 operators on $S^1\times B$ whose principal symbols are
 \[
 \sigma(A)|_{t=\pm\infty}(e^{i\theta},x,\tau,\xi)=\!\!\sum_{j+|\alpha|=N}a_{j,\alpha}^\pm(\theta,x)\tau^j\xi^\alpha,\ \  (e^{i\theta},x,\tau,\xi)\in S^*(S^1\times B).
 \]
 Having smooth symbols, their indices (and then $\delta_1([[A]_\ec]_1)$) can be more explicitly calculated as the integral of certain diferential forms
 on $S^*(S^1\times B)$.  

 With the aid of a computer algebra system, Patr\'icia Hess \cite{Hess} was able to follow this path and, using Fedosov's formula
 \cite[formula (1)]{FE}, find the index of a 2-by-2
 pseudodifferential operator on $S^1\times S^1$. That result was used in her proof that, for the case $B=S^1$, $K_0(\ac)\cong\Z^5$ and $K_1(\ac)\cong\Z^4$.


\end{document}